\documentclass[11pt,dvips,epsfig,reqno,graphics]{amsart}

\usepackage{a4wide,latexsym,amssymb,amsthm,amsmath,amsfonts,url}

\theoremstyle{plain}

\newtheorem{theorem}{Theorem}

\newtheorem{proposition}{Proposition}

\theoremstyle{remark}

\def\<{ \left < }
\def\>{ \right > }

\def\R{\mathbb{R}}

\def\E{\mathbb{E}}
\def\H{\mathbb{H}}
\def\S{\mathbb{S}}

\input epsfx.tex

\begin{document}

\title[Constant slope surfaces]
{From Golden Spirals to Constant Slope Surfaces\\
\ \\
\hfill{\scriptsize {\bf Motto:} EADEM MUTATA RESURGO}}

\author[M.~I.~Munteanu]{Marian Ioan Munteanu}

\thanks{This paper was supported by Grant PN-II ID 398/2007-2010 (Romania)}
\address{{Al.I.Cuza University of Iasi\\ Bd. Carol I, n. 11\\
700506 - Iasi\\ Romania\\
\newline\rm
\url{http://www.math.uaic.ro/~munteanu}}}
\email[]{{marian.ioan.munteanu@gmail.com}}

\begin{abstract}
In this paper, we find all constant slope surfaces in the Euclidean 3-space,
namely those surfaces for which the position vector of a point of the surface makes
constant angle with the normal at the surface in that point.
These surfaces could be thought as the bi-dimensional analogue of the generalized helices.
Some pictures are drawn by using the parametric equations we found.

\end{abstract}



\maketitle

\section{Introduction}

The study of spirals starts with the Ancient Greeks. A {\em logarithmic spiral} is a special type
of spiral discovered by Ren\'e Descartes and later extensively investigated by Jacob Bernoulli who called it
{\it spira mirabilis}. This curve has the property that the angle $\theta$ between its tangent and the
radial direction at every point is constant. This is the reason for which it is often known as
{\it equiangular spiral}. It can be described in polar coordinates $(\rho,\varphi)$ by the equation
$\rho=a{\rm e}^{b\varphi}$, where $a$ and $b$ are positive real constants with $b=\cot\theta$.
Of course, in the extreme cases, namely $\theta=\frac\pi2$ and $\theta\rightarrow0$,
the spiral becomes a circle of radius $a$, respectively it tends to a straight line.

\medskip

Why is this spiral miraculous? Due to its property that the size increases without altering the shape,
one can expect to find it on different objects around us. Indeed, let us mention here some phenomena where we may see curves
that are close to be a logarithmic spiral: the approach of a hawk to its prey, the approach of an
insect to a light source (see \cite{Boy99}), the arms of a spiral galaxy, the arms of the tropical cyclones,
the nerves of the cornea, several biological structures, e.g. {\it Romanesco broccoli}, {\it Convallaria majalis},
some spiral roses, sunflower heads, {\it Nautilus} shells and so on. For this reason these curves are named also
{\it growth spirals}. (See for other details \cite{Kap02, Muk04}.)

\medskip

A special kind of logarithmic spiral is the {\it golden spiral}, called often {\em a symbol a harmony and beauty}
due to its straight connection with Fibonacci numbers and certainly with the {\em golden ratio} $\Phi$.
It is clear that if we look around we see everyplace the {\it divine proportion}. Hence, it is not surprising the fact
that the golden spirals appear in many fields such as design (architecture, art, music, poetry), nature (plants, animals,
human, DNA, population growth), cosmology, fractal structures, markets, Theology and so on.

\smallskip

Let's take a look to its geometric property, namely the angle between the tangent direction in a point $p$ and the position vector
$\overrightarrow{p}$ is constant.

\smallskip

But how these curves appear in Nature?

For example, in Biology, after many observations, the trajectory of some primitive animals, in the way they orient
themselves by a pointlight source is assumed to be a spiral (see \cite{Boy99}). The mathematical model arising from
practice is proposed to be a curve making constant angle with some directions. The study of these trajectories can
be explained as follows.

First case, when we deal with a planar motion, using basic notions of the geometry of plane curves, it is proved that the parametric equations for
the trajectory of an insect flight are
$$
\left\{\begin{array}{l}
x(t)=ae^{t\cot\theta}\cos t\\[2mm]
y(t)=ae^{t\cot\theta}\sin t
\end{array}\right.
$$
which is an equiangular spiral.
(The light source is situated in the origin of an orthonormal frame and the starting point, at $t=0$, is $(a,0)$.)

Secondly, for a spatial motion, the condition that the radius vector and the tangent direction make a constant angle
yields an undetermined differential equation, so, further assumptions are needed. One of the supplementary
restrictions is for example that in its flight the insect keeps a constant angle between the direction to the light and the
direction up-down (see \cite{Boy99}). In this situation the curve one obtains is a {\em conchospiral} winding on a cone.

\smallskip

A similar technique is used in \cite{Lor06}, where two deterministic models for the flight of Peregrine Falcons and possible other raptors
as they approach their prey are examined by using tools in differential geometry of curves. Again, it is claimed that a certain angle
is constant, namely the angle of sight between falcon and prey. See also other papers related to this topic, e.g. \cite{Tuc01,TTAE}.

Another important problem comes from navigation: {\it find the} {\em loxodromes} {\it of the sphere}, namely those curves on $\S^2$
making a constant angle $\theta$ with the meridians.
Let us mention here the name of the Flemish geographer and cartographer Gerardus Mercator for his projection map
strictly related to loxodromes.
These curves are also known as {\em rhumb lines} (see \cite{Car76}).
If we take the parametrization of the sphere given by
$$
r(\phi,\psi)=(\sin\psi\cos\phi,\sin\psi\sin\phi,\cos\psi)
$$
with $\phi\in\left[0,2\pi\right)$ and $\psi\in\left[0,\pi\right]$, one gets the equations of loxodromes
$\log(\tan\frac\psi2)=\pm \phi\cot\theta$. Notice that $\phi=$constant represents a meridian, while
$\psi=$constant is the equation of a parallel. An interesting property of these curves is that they project, under
the stereographic projection, onto a logarithmic spiral.

A more general notion is the {\em slope line} also called {\em generalized helix} (see \cite{Kuh}). Recall that these curves
are characterized by the property that the tangent lines make constant angle with a fixed direction. Assuming that the torsion
$\tau(s)\neq0$, other two conditions characterize (independently) a generalized helix $\gamma$, namely
\begin{enumerate}
\item lines containing the principal normal $n(s)$ and passing through $\gamma(s)$ are parallel to a fixed plane,
or
\item $\frac{\kappa(s)}{\tau(s)}$ is constant, where $\kappa$ is the curvature of $\gamma$.
See also \cite{Car76}.
\end{enumerate}
\medskip

A lot of other interesting applications of helices are briefly enumerated in \cite{IB} (e.g. DNA double and collagen triple helix,
helical staircases, helical structure in fractal geometry and so on). All these make authors to say that the {\sl helix is one of the
most fascinated curves in Science and Nature.}

\smallskip

We pointed out so far that the curves making constant angle with the position vector are studied and have biological interpretation.
As a consequence, a natural question, more or less pure geometrical, is the following: {\sl Find all surfaces in the Euclidean space
making a constant angle with the position vector.} We like to believe that there exist such surfaces in nature due their spectacular forms.
We call these surfaces {\it constant slope surfaces}.

\smallskip

Concerning surfaces in the 3-dimensional Euclidean space making a constant angle with a fixed direction we notice that there exists a classification
of them. We mention two recent papers \cite{CS07} and \cite{MN}.
The applications of constant angle surfaces in the theory of liquid crystals and of layered fluids were
considered by P. Cermelli and A.J. Di Scala in \cite{CS07}, and they used for their study of surfaces
the Hamilton-Jacobi equation, correlating the surface and the direction field.
In \cite{MN} it is given another approach to classify all surfaces for which the unit normal makes a
constant angle with a fixed direction. Among them developable surfaces play an important role (see \cite{Nis09}).
Interesting applications of the usefulness of developable surfaces in architecture can be found in a very
attractive paper of G. Glaesner and F. Gruber in \cite{GG07}.

\smallskip

The study of constant angle surfaces was extended in different ambient spaces, e.g. in $\S^2\times\R$ (see \cite{DFVdVV}) and in $\H^2\times\R$
(see \cite{DM, DM1}). Here $\S^2$ is the unit 2-dimensional sphere and $\H^2$ is the hyperbolic plane. In higher dimensional Euclidean space,
hypersurfaces whose tangent space makes constant angle with a fixed direction are studied and a local description of
how these hypersurfaces are constructed is given. They are called helix hypersurfaces (see e.g. \cite{SR09}).

\smallskip

The main result of this paper is a classification theorem for all those surfaces for which the normal direction
in a point of the surface makes a constant angle with the position vector of that point. We find explicitly (in Theorem 1)
the parametric equations which characterize these surface. Roughly speaking, a constant slope surface can be constructed
by using an arbitrary curve on the sphere $\S^2$ and an equiangular spiral. More precisely, in any normal plane
to the curve consider a logarithmic spiral of a certain constant angle. One obtains a kind of fibre
bundle on the curve whose fibres are spirals.

The structure of this article is the following:
In {\em Preliminaries} we formulate the problem and we obtain the structure equations by using
formulas of Gauss and Weingarten. In Section 3 the embedding equations for constant slope surfaces
are obtained. Moreover some representations of such surfaces are given
in order to show their nice and interesting shapes.
In {\em Conclusion} we set the constant slope surfaces beside the other special surfaces.

\section{Preliminaries}

Let us consider $(M,g)\hookrightarrow (\tilde M, <~,~>)$ an isometric immersion of a manifold $M$
into a Riemannian manifold $\tilde M$ with Levi Civita connection ${\stackrel{\circ}{\nabla}}$. Recall
the formulas of Gauss and Weingarten

\qquad
{\bf (G)} ${\stackrel{\circ}{\nabla}}_XY=\nabla_XY+{\sf{h}}(X,Y)$

\qquad
{\bf (W)} ${\stackrel{\circ}{\nabla}}_XN=-A_NX+\nabla^\bot_XN$

for every $X$ and $Y$ tangent to $M$ and for every $N$ normal to $M$. Here $\nabla$ is
the Levi Civita connection on $M$, ${\sf h}$ is a symmetric $(1,2)$-tensor field taking values
in the normal bundle and called {\it the second fundamental form} of $M$,
$A_N$ is the shape operator associated to $N$ also known as the Weingarten operator corresponding to $N$
and $\nabla^\bot$ is the induced connection in the normal bundle of $M$.
We have $<{\sf h}(X,Y),N>=g(A_NX,Y)$, for al $X,Y$ tangent to $M$, and hence $A_N$ is a
symmetric tensor field.

Now we consider an orientable surface $S$ in the Euclidean space $\E^3\setminus\{0\}$.
For a generic point $p$ in $\E^3$ let us use the same notation, namely $p$, also for
its position vector.
We study those surfaces $S$ in $\E^3\setminus\{0\}$ making a constant angle $\theta$ with $p$.

Denote by $<~,~>$ the Euclidean metric on $\E^3\setminus\{0\}$, by ${\stackrel{\circ}{\nabla}}$ its flat connection
and let $g$ be the restriction of $<~,~>$ to $S$.

Let $\theta\in[0,\pi)$ be a constant.

Two particular cases, namely $\theta=0$ and $\theta=\frac\pi2$ will be treated separately.

Let $N$ be the unit normal to the surface $S$ and $U$ be the projection of $p$ on the
tangent plane (in $p$) at $S$. Then, we can decompose $p$ in the form
\begin{equation}
\label{eq:1}
p=U+a N, \quad a\in C^\infty(S).
\end{equation}
Denoting by $\mu=|p|$ the length of $p$, one gets $a=\mu\cos\theta$ and hence
$p=U+\mu\cos\theta N$. It follows
\begin{equation}
\label{eq:2}
|U|=\mu\sin\theta.
\end{equation}
Since $U\neq 0$ one considers $e_1=\frac U{|U|}$, unitary and tangent to $S$ in $p$.
Let $e_2$ be a unitary tangent vector (in $p$) to $S$ and orthogonal to $e_1$.
We are able now to write
\begin{equation}
\label{eq:3}
p=\mu(\sin\theta e_1 + \cos\theta N).
\end{equation}

For an arbitrary vector field $X$ in $\E^3$ we have
\begin{equation}
\label{eq:4}
{\stackrel{\circ}{\nabla}}_Xp=X.
\end{equation}

If $X$ is taken to be tangent to $S$, a derivation in \eqref{eq:3} with respect to $X$, combined with the relation above yield
\begin{equation}
\label{eq:5}
X  =  \sin\theta\left(X(\mu)e_1+\mu\nabla_Xe_1+\mu h(X,e_1)N\right)
 + \cos\theta \left(X(\mu)N-\mu AX\right)
\end{equation}
where $h$ is the real valued second fundamental form (${\sf h}(X,Y)=h(X,Y)N$) and $A=A_N$ is the Weingarten operator on $S$.

Identifying the tangent and the normal parts respectively we obtain
\begin{equation}
\label{eq:6}
X=\sin\theta\big(X(\mu)~e_1+\mu\nabla_Xe_1\big) - \mu\cos\theta AX
\end{equation}
\begin{equation}
\label{eq:7}
0= \mu \sin\theta h(X,e_1) + \cos\theta X(\mu).
\end{equation}

We plan to write the expression of the Weingarten operator in terms of the local basis
$\{e_1,e_2\}$. Since ${\rm grad} ~\mu =\frac 1\mu ~p$ one obtains
\begin{equation}
\label{eq:9}
X(\mu)=<X,e_1>\sin\theta.
\end{equation}

Thus, the relation \eqref{eq:7} becomes
$$
0=\mu\sin\theta g(Ae_1,X)+\cos\theta\sin\theta g(e_1,X),\quad \forall X\in\chi(S)
$$
and hence one gets
\begin{equation}
\label{eq:11}
Ae_1=-\frac{\cos\theta}\mu~e_1
\end{equation}
namely $e_1$ is a principal direction for the Weingarten operator.

\medskip

Consequently, there exists also a smooth function $\lambda$ on $S$ such that
\begin{equation}
\label{eq:14}
Ae_2=\lambda e_2.
\end{equation}
Hence, the second fundamental form $h$ can be written as
$\left(
\begin{array}{cc}
-\frac{\cos\theta}\mu & 0\\
0 & \lambda
\end{array}
\right)$.

\begin{proposition}
\label{prop:1}
The Levi Civita connection on $S$ is given by
\begin{equation}
\label{eq:LC}
\nabla_{e_1}e_1=0 ~,~
\nabla_{e_1}e_2=0 ~,~
\nabla_{e_2}e_1=\frac {1+\mu\lambda\cos\theta}{\mu\sin\theta} ~ e_2 ~,~
\nabla_{e_2}e_2=-\frac{1+\mu\lambda\cos\theta}{\mu\sin\theta} ~ e_1.
\end{equation}
\end{proposition}
\proof
Considering $X=e_1$ (respectively $X=e_2$) in \eqref{eq:6} and \eqref{eq:9}, and combining with \eqref{eq:11} one obtains
\eqref{eq:LC}${}_1$ (respectively \eqref{eq:LC}${}_3$). The other two statements can be obtained immediately.

\endproof

\section{The characterization of constant slope surfaces}

Let $r:S\longrightarrow \E^3\setminus\{0\}$ be an isometric immersion of the surface $S$ in the
Euclidean space from which we extract the origin.

As a consequence of the Proposition \ref{prop:1} we have
$[e_1,e_2]=-\frac {1+\mu\lambda\cos\theta}{\mu\sin\theta} ~ e_2$.
Hence, one can consider local coordinates $u$ and $v$ on $S$ such that
$\partial_u\equiv\frac\partial{\partial u}=e_1$ and $\partial_v\equiv\frac\partial{\partial v}=\beta(u,v) e_2$
with $\beta$ a smooth function on $S$.
It follows that $\beta$ should satisfy
\begin{equation}
\label{eq:20}
\beta_u-\beta\frac {1+\mu\lambda\cos\theta}{\mu\sin\theta}=0.
\end{equation}

The Schwarz equality $\partial_u\partial_v N=\partial_v\partial_u N$ yields
$$
(\beta\lambda)_u e_2=\frac{\mu_v}{\mu^2} \cos\theta e_1 -
       \beta \frac{\cos\theta(1+\mu\lambda\cos\theta)}{\mu^2\sin\theta} e_2.
$$
Consequently $\mu=\mu(u)$ and
\begin{equation}
\label{eq:23}
\mu^2\sin\theta\lambda_u+(1+\mu\lambda\cos\theta)(\mu\lambda+\cos\theta)=0.
\end{equation}

Moreover, from $\mu(u)^2=|r(u,v)|^2$ taking the derivative with respect to $u$ and using that
$r_u=e_1$ and that the tangent part of the position vector $r(u,v)$ is equal to $\mu(u)\sin\theta e_1$, we obtain
$\mu'(u)=\sin\theta$. It follows that $\mu(u)=u\sin\theta+\mu_0$, where $\mu_0$ is a real constant
which can be taken zero after a translation of the parameter $u$. So, we can write
\begin{equation}
\label{eq:24}
\mu(u)=u\sin\theta.
\end{equation}
Accordingly, the PDE \eqref{eq:23} becomes
\begin{equation}
u^2\sin^3\theta\lambda_u+(1+\lambda u\sin\theta\cos\theta)(\lambda u\sin\theta+\cos\theta)=0.
\end{equation}
Denoting $\rho(u,v)=\lambda(u,v)u\sin\theta$, the previous equation comes out
$$
   u\rho_u\sin^2\theta=-\cos\theta \big(\rho^2+2\rho\cos\theta+1\big).
$$
By integration one gets
\begin{equation}
\label{eq:S1}
\rho(u,v)=-\cos\theta-\sin\theta\tan\big(\cot\theta\log u+Q(v)\big)
\end{equation}
where $Q$ is a smooth function depending on the parameter $v$.
From here we directly have $\lambda$. In the sequel, we
solve the differential equation \eqref{eq:20} obtaining
\begin{equation}
\label{eq:S2}
\beta(u,v)=u\cos\big(\cot\theta\log u+Q(v)\big)\varphi(v)
\end{equation}
where $\varphi$ is a smooth function on $S$ depending on $v$.

We can state now the main result of this paper.
\begin{theorem}
\label{th:1}
Let $r:S\longrightarrow \R^3$ be an isometric immersion of a surface $S$ in the Euclidean $3$-dimensional space.
Then $S$ is of constant slope if and only if either it is an open part of the Euclidean $2$-sphere centered in the origin,
or it can be parametrized by
\begin{equation}
\label{eq:imm}
r(u,v)=u\sin\theta \big(\cos\xi ~ f(v) + \sin\xi ~ f(v)\times f'(v) \big)
\end{equation}
where $\theta$ is a constant (angle) different from $0$, $\xi=\xi(u)=\cot\theta \log u$ and $f$ is a unit speed curve on the Euclidean sphere $\S^2$.
\end{theorem}
\proof
(Sufficiency)
As regards the sphere, it is clear that the position vector is normal to the surface.
The angle is zero in this case.

Concerning the parametrization \eqref{eq:imm}, let us prove now that it characterizes a constant slope surface in $\R^3$.

First we have to compute $r_u$ and $r_v$ in order to determine the tangent plane to $S$.
Hence

$r_u=-\sin(\xi-\theta) ~ f(v) + \cos(\xi-\theta) ~ f(v)\times f'(v)$

$r_v=u\sin\theta(\cos\xi+\kappa_g\sin\theta)~f'(v)$,

where $\kappa_g$ is the proportionality factor between $f\times f''$ and $f'$ and it
represents, up to sign, the geodesic curvature of the curve $f$.

The unit normal to $S$ is
$$
N=\cos(\xi-\theta)~f(v) + \sin(\xi-\theta) ~ f(v)\times f'(v).
$$
Then, the angle between $N$ and the position $r(u,v)$ is given by
$$
\cos\widehat{(N,r)}=\cos\theta=constant.
$$
Thus the sufficiency is proved.

\smallskip

Conversely, we start with the local coordinates $u$ and $v$, the metric $g$ on $S$ and its Levi Civita connection $\nabla$.
In order to determine the immersion $r$, we have to exploit the formula of Gauss.
More precisely, we have first
$r_{uu}=\nabla_{e_1}e_1+h(e_1,e_1) N$
and hence
$$
r_{uu}=-\frac{\cot\theta}u~N.
$$
Then we know that the position vector decomposes as
$$
r(u,v)=u\sin^2\theta r_u + u\sin\theta\cos\theta N.
$$
These two relations yield the following differential equation
\begin{equation}
\label{eq:P15}
r-u\sin^2\theta r_u+u^2\sin^2\theta r_{uu}=0.
\end{equation}
It follows that there exist two vectors $\tilde c_1(v)$ and $\tilde c_2(v)$ such that
\begin{equation}
\label{eq:P3}
r(u,v)=u\cos\big(\cot\theta \log u\big) \tilde c_1(v)+
          u\sin\big(\cot\theta \log u\big) \tilde c_2(v)
\end{equation}
Denote, for the sake of simplicity, by $\xi=\xi(u)=\cot\theta \log u$.

Computing $|r(u,v)|=\mu(u)=u\sin\theta$ one obtains
$$
\sin^2\theta=\frac{1+\cos 2\xi}2 ~ |\tilde c_1(v)|^2+\frac{1-\cos 2\xi}2 ~ |\tilde c_2(v)|^2
             +\sin 2\xi <\tilde c_1(v),\tilde c_2(v)>
$$
for any $u$ and $v$.
Accordingly, $|\tilde c_1(v)|=|\tilde c_2(v)|$, $|\tilde c_1(v)|+|\tilde c_2(v)|=2\sin^2\theta$,
$<\tilde c_1(v),\tilde c_2(v)>=0$.
As consequence there exist two unitary and orthogonal vectors, $c_1(v)$ and $c_2(v)$,
such that $\tilde c_1(v)=\sin\theta c_1(v)$ and $\tilde c_2(v)=\sin\theta c_2(v)$. At this point
the condition $|r_u|^2=1$ is satisfied.

The immersion $r$ becomes
\begin{equation}
\label{eq:P6p}
r(u,v)=u\sin\theta \big(\cos\xi c_1(v) + \sin\xi c_2(v)\big).
\end{equation}

From the orthogonality of $r_u$ and $r_v$, one gets $<c_1'(v),c_2(v)>=0$, for any $v$.
Thus, due the fact that $c_2$ is orthogonal both to $c_1$ and to $c_1'$, it results to be collinear with
the cross product $c_1 \times c_1'$. Moreover we have
$$
c_2'(v)=-\frac{(c_1,c_1',c_1'')}{|c_1'|^3} ~ c_1'(v).
$$
A straightforward computation shows us that $\frac{(c_1,c_1',c_1'')}{|c_1'|^3}=\tan Q(v)$
and $\varphi(v)=\frac{\sin\theta |c_1'(v)|}{\cos Q(v)}$
(see relations \eqref{eq:S1}, \eqref{eq:S2}).

In the sequel, the equation
$r_{uv}=\beta_u e_2 + \nabla_{e_1}e_2 + h(e_1,e_2) N$
furnishes no more information.

Let us make a change of parameter $v$, namely suppose that $|c_1'(v)|=1$. Therefore $c_1$ can
be thought of as a unit curve on the sphere $\S^2$ and denote it by $f$. Consequently, $c_2=f\times f'$.
With these notations we have the statement of the theorem.

Let us notice that the last equation, namely
$r_{vv}=-\beta^2\frac{1+\lambda u\sin\theta\cos\theta}{u\sin^2\theta} ~ r_u +
    \frac{\beta_v}\beta ~ r_v + \beta^2\lambda N$
gives no other conditions.

It is interesting to mention separately the case $\theta =\frac\pi2$, even that it derives from the
general parametric equation of the surface. If this happens, the surface $S$ can be parametrized as
$$r(u,v)=u f(v)$$
which represents
\begin{itemize}
\item [ (i)] either the equation of a cone with the vertex in the origin,
\item[(ii)] or the equation of an open part of a plane passing through the origin.
\end{itemize}
Here $f$ represents, as before, a unit speed curve on the Euclidean 2-sphere.

To end the proof the particular case $\theta=0$ has to be analyzed.

Case $\theta=0$. Then the position vector is normal to the surface. We have
$<r,r_u>=0$ and $<r,r_v>=0$. It follows that $|r|^2=constant$, i.e. $S$ is
an open part of the Euclidean $2$-sphere.

\endproof

Let us remark that the only flat constant slope surfaces are open parts of planes and cones,
case in which the angle between the normal and the position vector is a right angle.
Similarly the only minimal constant slope surfaces in $\E^3$ are open parts of planes.

\medskip

As we expected, the shape of these surfaces is quite interesting and we give some pictures made
with Matlab in order to have an idea of what they look like.

\begin{figure}[htb]
\begin{center}
\epsfxsize=90mm \centerline{\leavevmode \epsffile{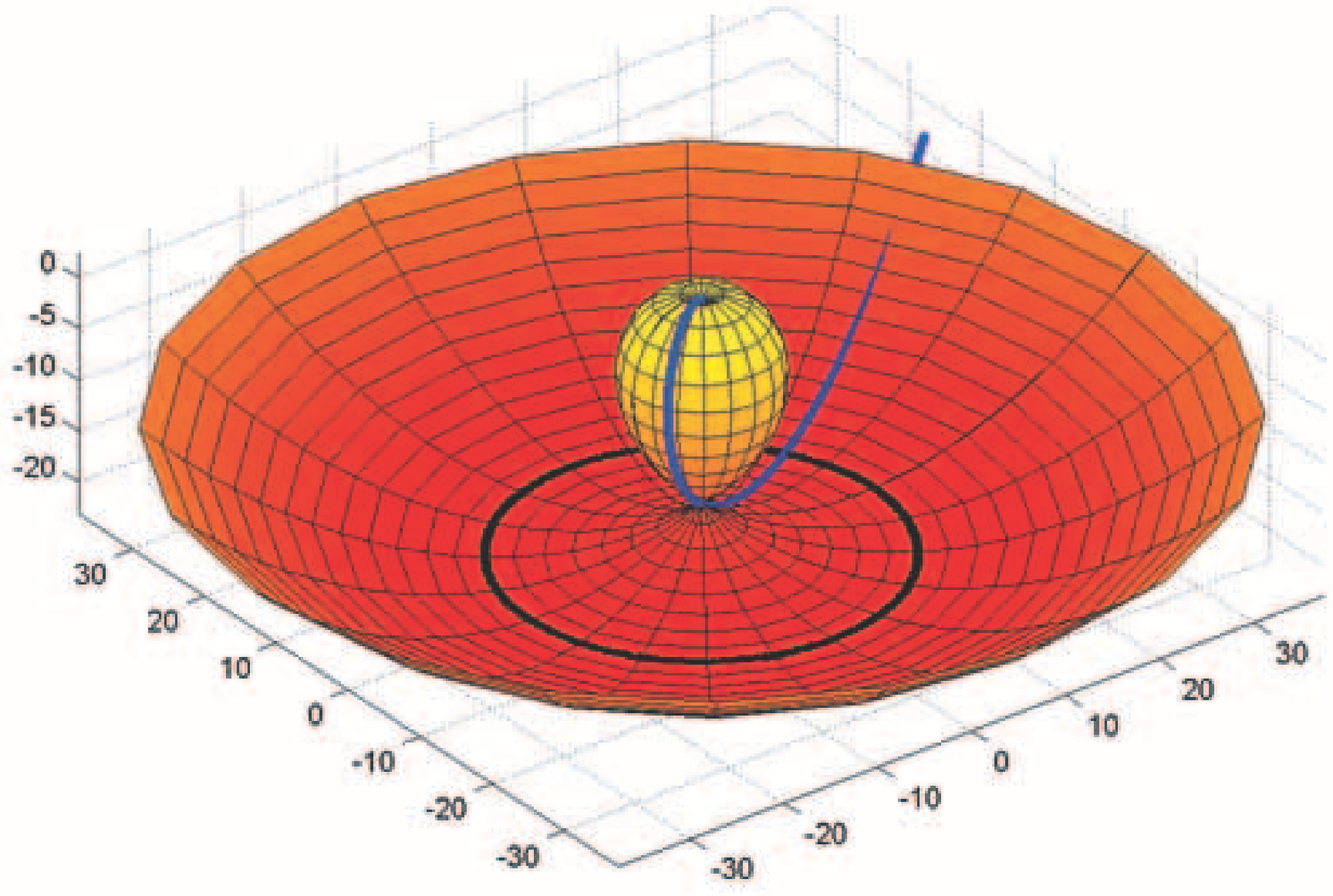}
\epsfxsize=90mm \epsffile{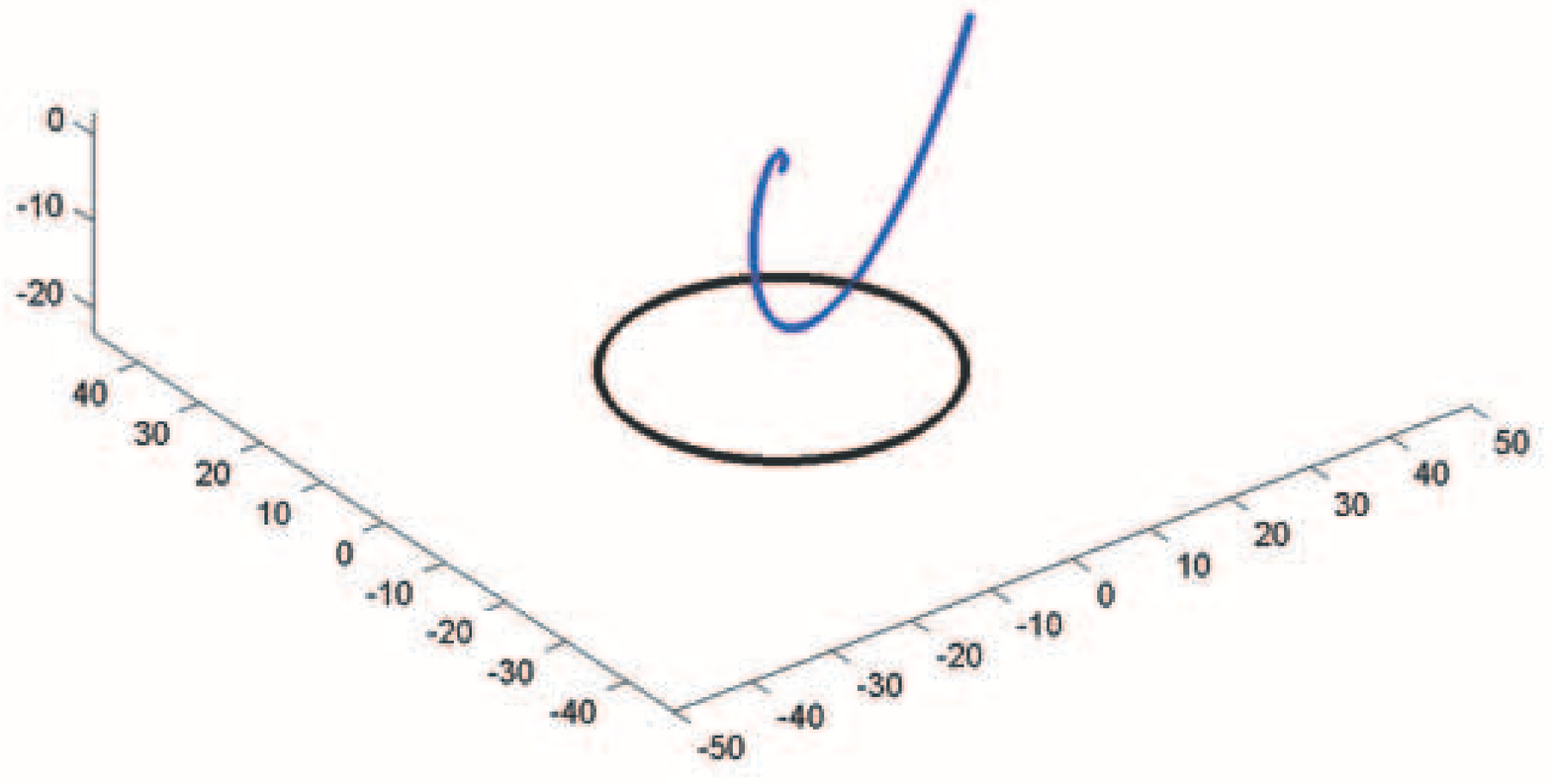} }
\end{center}
\*{$\theta=\frac\pi5$\hspace{10mm} $f(v)=(\cos v,\sin v,0)$}
\end{figure}

Let us point more attention to this picture (but not necessary with $\theta=\frac\pi5$), when $f(v)=(\cos v,\sin v,0)$. Then $f(v)\times f'(v)=(0,0,1)$ for all $v$ and consequently
the slope surface is parametrized by
$$r(u,v)=u\sin\theta\ (\cos(\xi(u))\cos v,\cos(\xi(u))\sin v,\sin(\xi(u))).$$
It is interesting to notice that the parametric line $v=0$ is an equiangular spiral in the $(xz)$-plane
and it is drawn in blue. The black line represents the parametric line $u=u_0$ (for a given $u_0$).

\pagebreak

In the following we represent other two slope surfaces ($\theta\neq\frac\pi2$) together with their parametric lines in a certain point.


\begin{figure}[htb]
\begin{center}
\epsfxsize=90mm \centerline{\leavevmode \epsffile{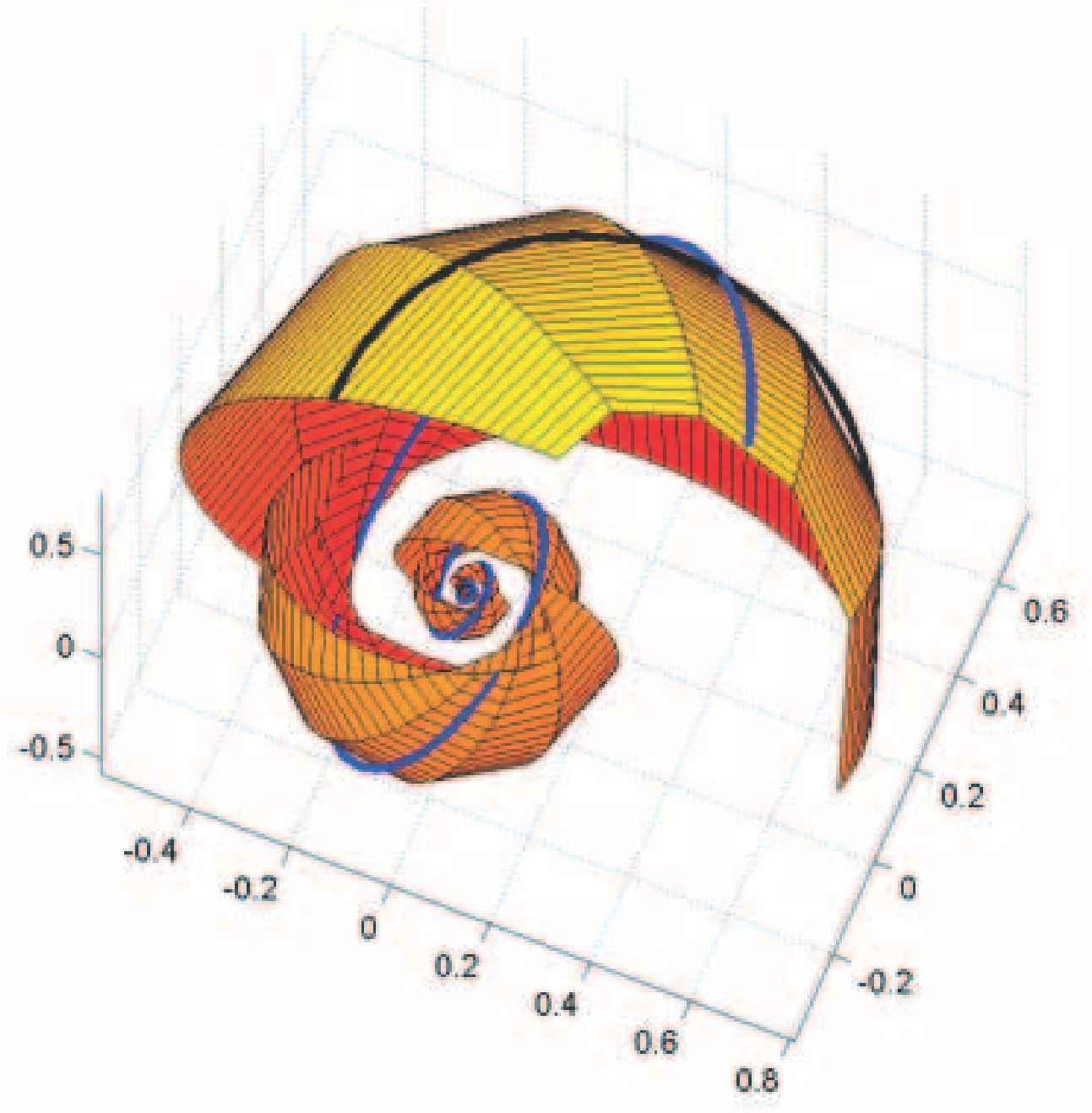}
\epsfxsize=90mm \epsffile{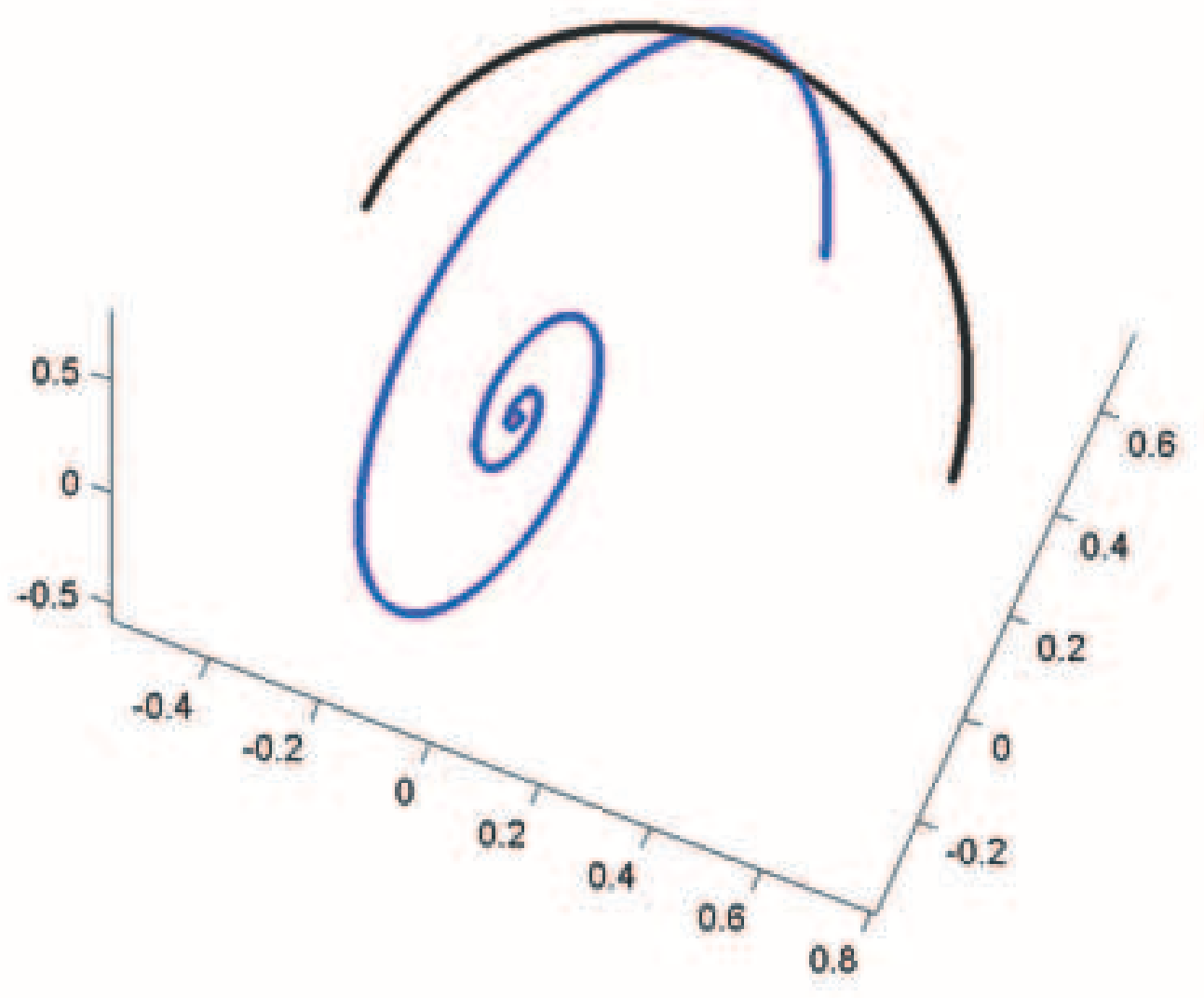} }
\end{center}
\*{$\theta=\frac\pi{15}$\hspace{5mm} $f(v)=(\cos^2v,\cos v\sin v,\sin v)$}
\end{figure}


\begin{figure}[htb]
\begin{center}
\epsfxsize=90mm \centerline{\leavevmode \epsffile{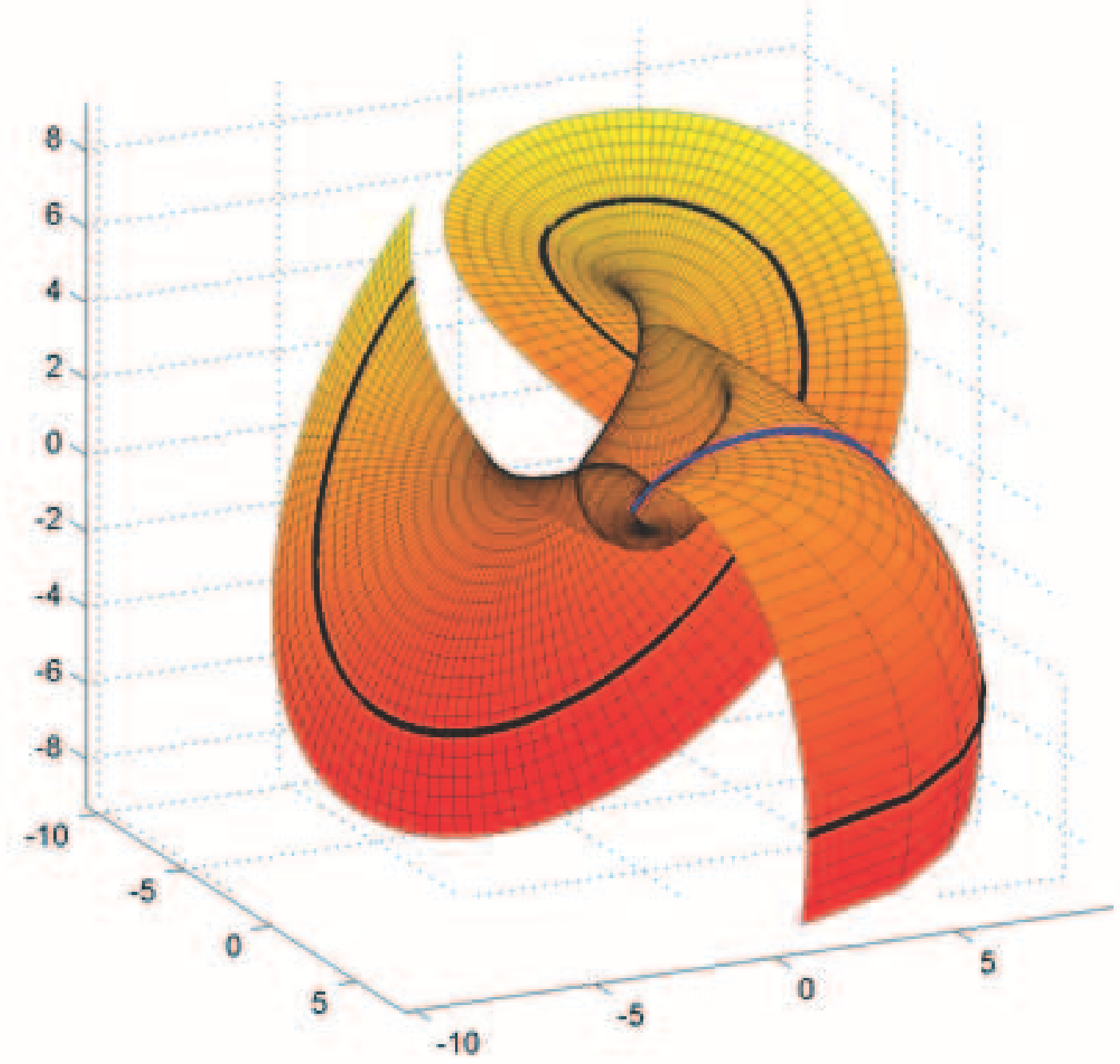}
\epsfxsize=90mm \epsffile{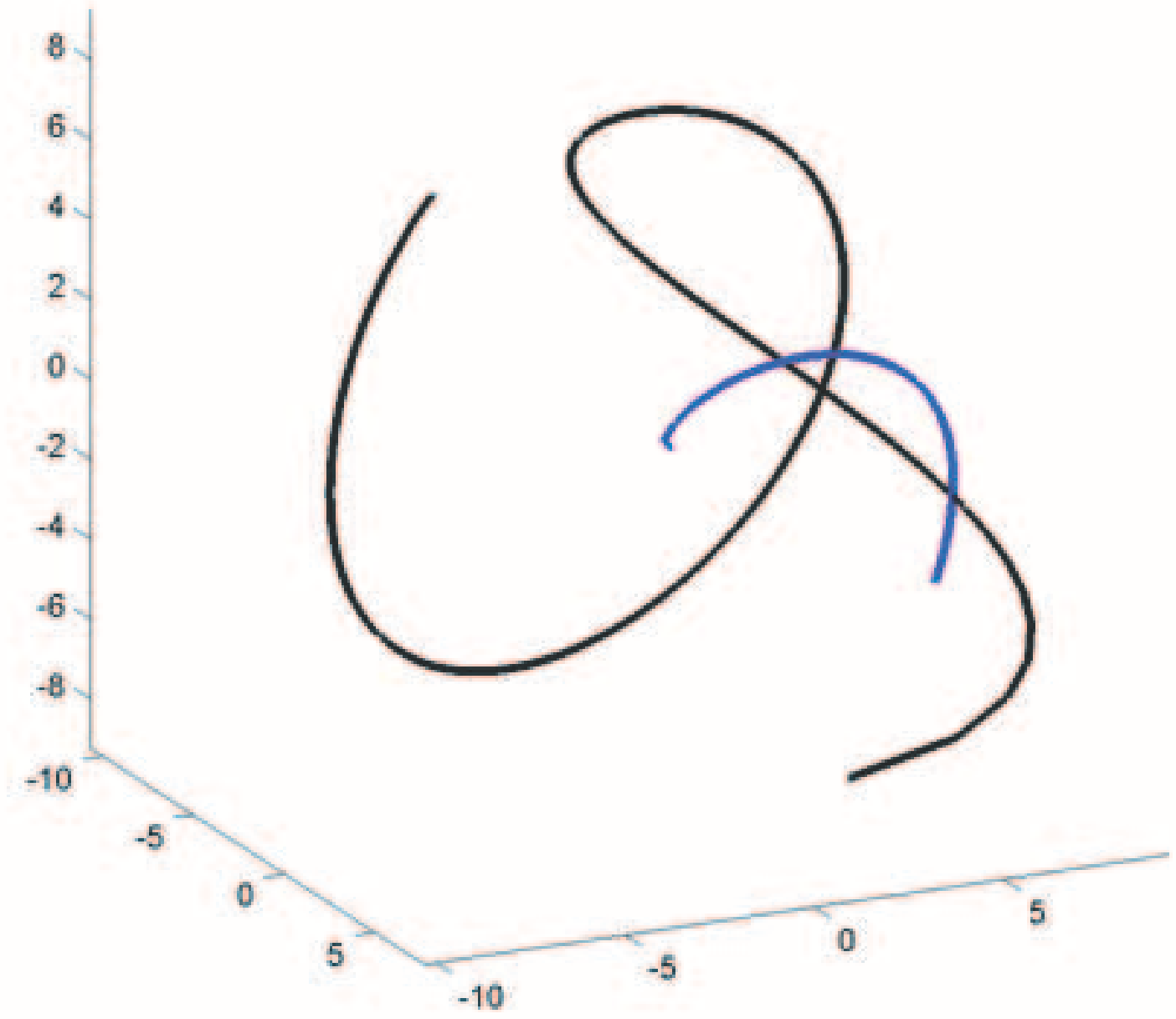} }
\end{center}
\*{$\theta=\frac\pi4$ \hspace{5mm} $f(v)=(\sin(\psi(v))\cos(\phi(v)),\sin(\psi(v))\sin(\phi(v)),\cos(\psi(v)))$}\\
\*{$\phi(v)=\log(v)$ \hspace{5mm} $\psi(v)=v$}
\end{figure}

\pagebreak

And now, here we are the cone, a slope surface with the constant angle $\frac\pi2$. One can see the spiral, as parametric line, degenerated into a straight line
on the cone.


\begin{figure}[htb]
\begin{center}
\epsfxsize=90mm \centerline{\leavevmode \epsffile{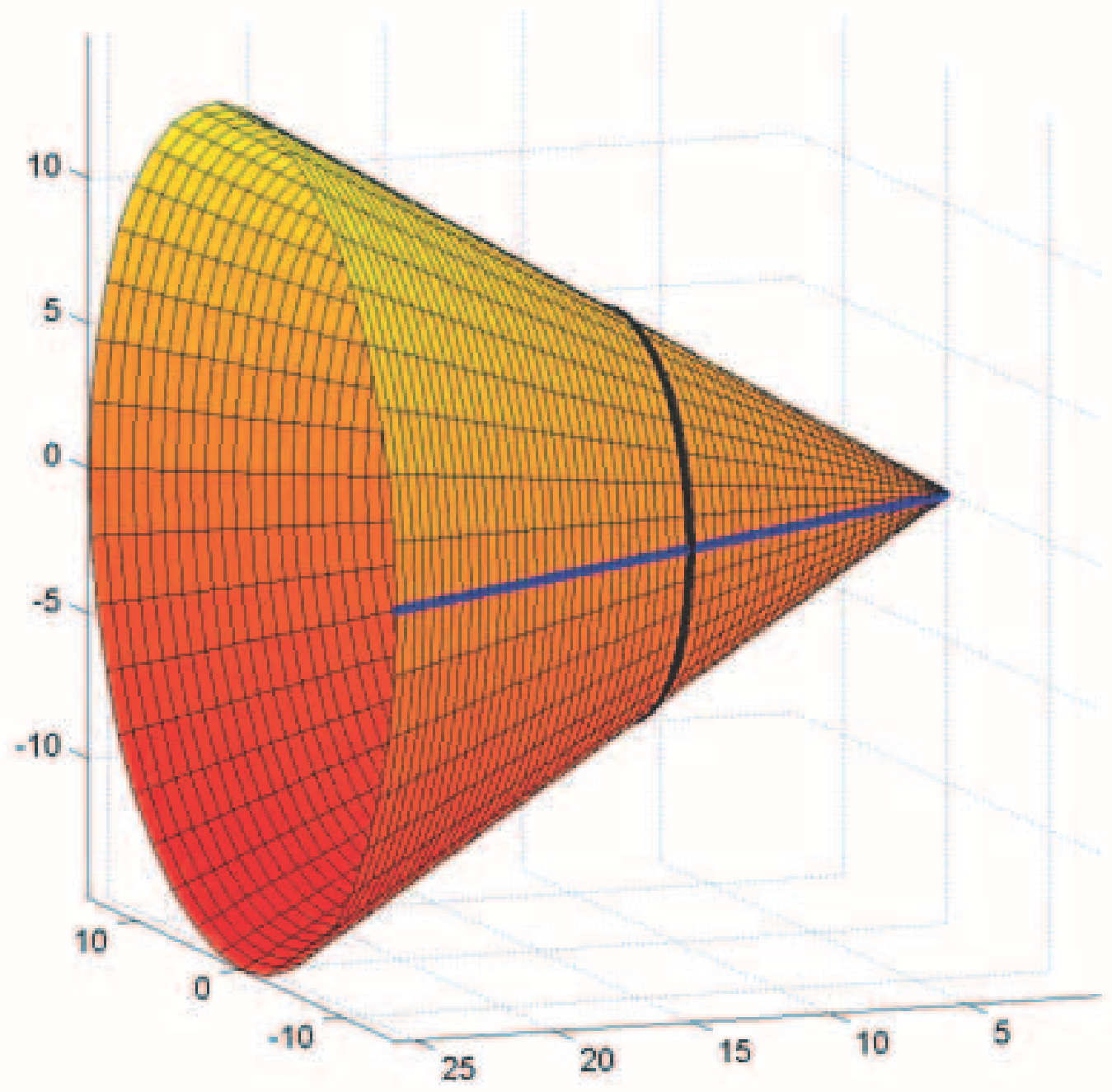}
\epsfxsize=90mm \epsffile{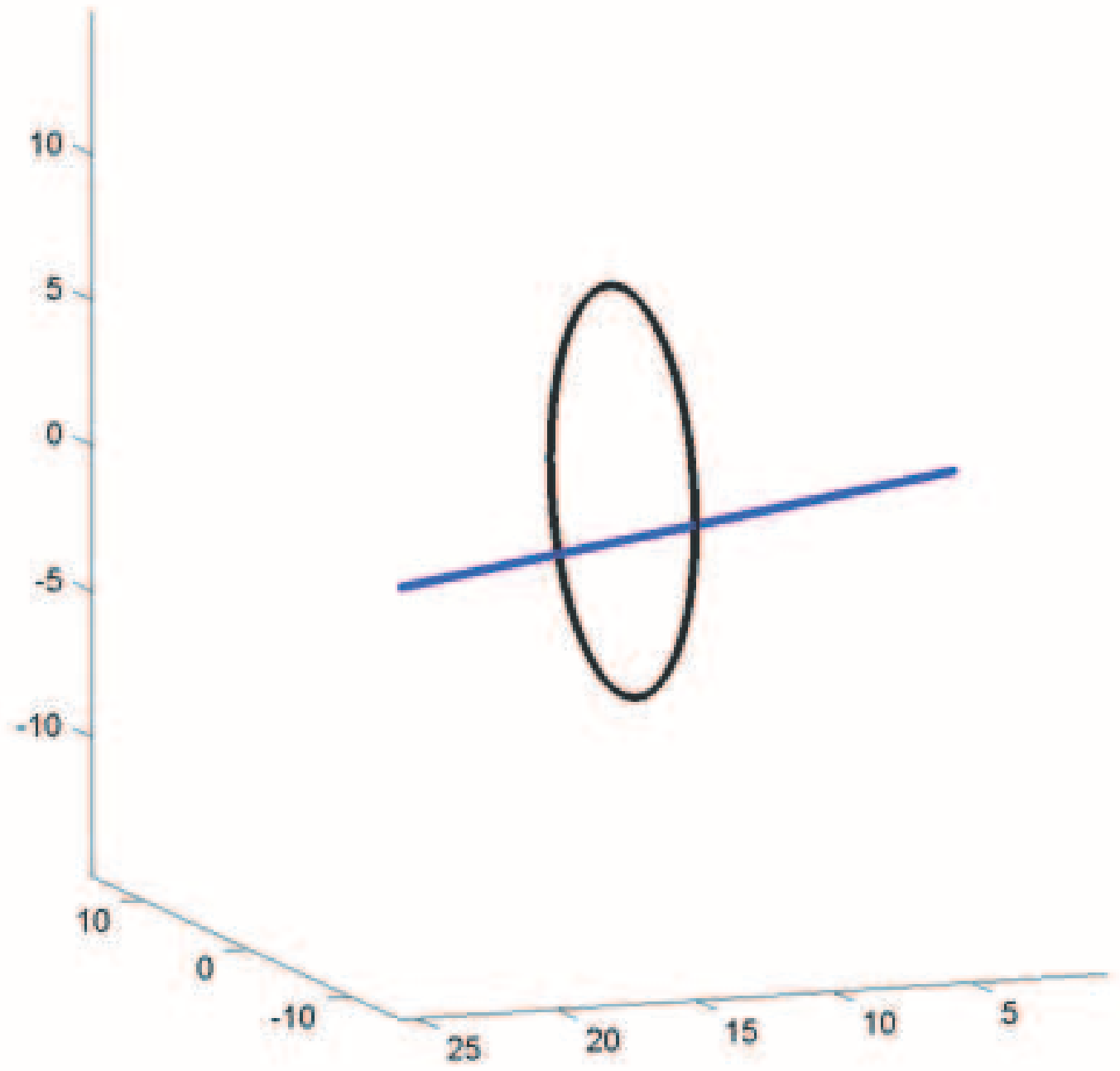} }
\end{center}
\*{$\theta=\frac\pi2$ \hspace{5mm} $f(v)=\frac12\ (\cos v,\sqrt{3},\sin v)$}
\end{figure}


\section{Conclusion}

Surfaces for which the normal in a point makes constant angle with the position vector have nice shapes and they are
interesting from the geometric point of view.
The study of these surfaces is similar, on one hand, to that of the logarithmic spirals and generalized helices
(also known as slope lines) and on the other hand, it is close related to the study of constant angle surfaces and spiral surfaces.
These last ones were introduced at the end of the nineteenth century by Maurice L\'evy \cite{Lev878}.
At least for their shapes, one can say that constant slope surfaces are one of the most fascinated surfaces in the Euclidean 3-space.


{\small
{\bf  Acknowledgements.} I would like to thank Prof. Antonio J. di Scala for reading very carefully an earlier draft
and for his valuable advices and suggestions concerning the subject of this paper.}

\end{document}